\documentclass[12pt]{article}
\usepackage{amssymb}
\usepackage{epsfig}
\usepackage{color}

\normalbaselines

\newtheorem{Theorem}{Theorem}[section]

\newtheorem{Proposition}[Theorem]{Proposition}

\newtheorem{Lemma}[Theorem]{Lemma}

\newtheorem{Corollary}[Theorem]{Corollary}

\newtheorem{Remark}[Theorem]{Remark}

\newcommand{\RR}{{{\rm I} \kern -.15em {\rm R} }}

\newcommand{\C}{{{\rm l} \kern -.42em {\rm C} }}

\newcommand{\nat}{{{\rm I} \kern -.15em {\rm N} }}

\newcommand{\be}{\begin{equation}}
\newcommand{\ee}{\end{equation}}
\newcommand{\beq}{\begin{eqnarray}}
\newcommand{\eeq}{\end{eqnarray}}
\newcommand{\beqs}{\begin{eqnarray*}}
\newcommand{\eeqs}{\end{eqnarray*}}
\newcommand{\bt}{\begin{Theorem}}
\newcommand{\et}{\end{Theorem}}
\newcommand{\br}{\begin{Remark}}
\newcommand{\er}{\end{Remark}}
\newcommand{\bc}{\begin{Corollary}}
\newcommand{\ec}{\end{Corollary}}
\newcommand{\bl}{\begin{Lemma}}
\newcommand{\el}{\end{Lemma}}
\newcommand{\bd}{\begin{definition}}
\newcommand{\ed}{\end{definition}}

\title{Exponential stability of the wave equation with\\ memory
and time delay}
\author{{\sc F. Alabau-Boussouira}\\
{\small LMAM, Universit\'{e} de Lorraine and CNRS (UMR 7122)}\\
{\small 57045 Metz Cedex 1, France}
\\\\
{\sc S. Nicaise}\\
{\small LAMAV, FR CNRS 2956, Institut des Sciences et Techniques 
de Valenciennes}\\
{\small 
Universit\'{e} de Valenciennes et du Hainaut Cambr\'{e}sis}\\
{\small 59313 Valenciennes Cedex 9, France}\\\\
{\sc C. Pignotti}
\\
{\small 
Dipartimento di Ingegneria e Scienze dell'Informazione e Matematica}\\
{\small Universit\`{a} di L'Aquila,
67010 L'Aquila, Italy}}

\date{}

\begin{document}

\textwidth=160 mm

\textheight=225mm

\parindent=8mm

\frenchspacing

\maketitle

\bigskip

\begin{abstract}
We study the asymptotic behaviour of the wave equation with viscoelastic damping
in presence of a 
time--delayed damping.
We prove exponential stability if the amplitude of the time delay term is small enough.
\end{abstract}

\vspace{5 mm}

\def\qed{\hbox{\hskip 6pt\vrule width6pt
height7pt
depth1pt  \hskip1pt}\bigskip}



\section{Introduction}
\label{intro}\hspace{5mm}

\setcounter{equation}{0}

This paper is devoted to the stability analysis of a viscoelastic model.
In particular, we consider a model combining viscoelastic damping and time-delayed damping. We prove an exponential stability result provided that the amplitude of time-delayed damping is small enough.  Moreover, we give a precise 
estimate on this smallness condition. This shows that even if delay effects usually generate instability (see e.g. \cite{Datko, DLP, NPSicon06, XYL}), the damping due to viscoelasticity can counterbalance them.

Let $\Omega\subset\RR^n$ be an open bounded set with a smooth boundary.
Let us consider the following problem:

\begin{eqnarray}
& &u_{tt}(x,t) -\Delta u (x,t)-
\int_0^\infty \mu (s)\Delta u(x,t-s) ds\nonumber \\
& &\hspace{5 cm}
 +
k u_t(x,t-\tau )=0\quad \mbox{\rm in}\ \Omega\times
(0,+\infty)\label{1.1d}\\
& &u (x,t) =0\quad \mbox{\rm on}\ \partial\Omega\times
(0,+\infty)\label{1.2d}\\
& &u(x,t)=u_0(x, t)\quad \hbox{\rm
in}\ \Omega\times (-\infty, 0] \label{1.3d}
\end{eqnarray}
where  the initial datum $u_0$ belongs to a suitable space, the constant $\tau >0$ is the time delay, $k$ is a real number and
the  memory kernel $\mu :[0,+\infty)\rightarrow [0,+\infty)$ 
is a locally absolutely continuous function
satisfying

i) $\mu (0)=\mu_0>0;$

ii) $\int_0^{+\infty} \mu (t) dt=\tilde \mu <1;$

iii) $\mu^{\prime} (t)\le -\alpha \mu (t), \quad \mbox{for some}\ \ \alpha >0.$

\noindent
We know that the above problem is exponentially stable for $k =0$ (see e.g.
\cite{GiorgiRiveraPata}).

We will show that an exponential stability result holds if the delay parameter $k$ is {\em small} with  respect to the memory kernel.

Observe that for $\tau =0$ and $k>0$ the   model $(\ref{1.1d})-(\ref{1.3d})$
presents both viscoelastic and standard dissipative damping. Therefore, in that case, under the above assumptions on  the kernel 
$\mu ,$ the model is exponentially stable. 

We will see that 
exponential stability also occurs for $k<0,$ under a suitable smallness assumption on  $\vert k\vert.$  Note that the term $ku_t(t)$ with $k<0$ is  a so--called
anti--damping (see e.g. \cite{FZ}), namely a damping with an opposite sign
with respect to the standard dissipative one, and therefore it induces instability.
Indeed, in absence of viscoelastic damping, i.e. for $\mu\equiv 0,$ the solutions of the above problem, with $\tau =0$ and $k<0,$ grow exponentially to infinity.

We will prove our stability results by using a perturbative approach, first introduced in \cite{SCL12} (see also \cite{NiPi2013} for a more general setting).

The stabilization problem for model $(\ref{1.1d})-(\ref{1.3d})$ has been studied
also by Guesmia in \cite{Guesmia} by using a different approach based on the construction of a suitable Lyapunov functional.
Our analysis allows to determine an explicit estimate on the constant $k_0$ (cf. Theorem \ref{stab2}).
Moreover, our approach can be extended to the case of localized viscoelastic
damping (cf. \cite{Rivera}). In fact, we first prove the exponential stability
of an auxiliary problem having a decreasing energy and then, regarding the original problem as a perturbation of that one,  we extend
the exponential decay estimate to it.

The paper is organized as follows. In sect. 2 we study the well--posedness
by introducing an appropriate functional setting and we formulate our stability result. In sect. 3 we introduce the auxiliary problem and prove the exponential decay estimate for it. Then, the stability result is extended to the original problem.  
\eject
\section{Main results and preliminaries}

\label{pbform}\hspace{5mm}

\setcounter{equation}{0}

As in \cite{Dafermos}, let us introduce the new variable
\begin{equation}\label{eta}
\eta^t(x,s):=u(x,t)-u(x,t-s).
\end{equation}

\noindent
Moreover, 
as in \cite{NPSicon06}, we define
\begin{equation}\label{zeta}
z(x,\rho ,t):=
u_t(x,t-\tau\rho),\quad x\in\Omega,\ \rho\in (0,1),\ t>0.
\end{equation}

\noindent
Using (\ref{eta}) and (\ref{zeta}) we can rewrite (\ref{1.1d})--(\ref{1.3d})
as

\begin{eqnarray}
& &u_{tt}(x,t)= (1-\tilde \mu)\Delta u (x,t)+
\int_0^\infty \mu (s)\Delta \eta^t(x,s) ds\nonumber\\
& &\hspace{5 cm}
 -k
z(x,1,t)\quad \mbox{\rm in}\ \Omega\times
(0,+\infty)\label{e1d}\\
& & \eta_t^t(x,s)=-\eta^t_s(x,s)+u_t(x,t)\quad \mbox{\rm in}\ \Omega\times
(0,+\infty)\times (0,+\infty ),\label{e2d}\\
& &\tau z_t(x,\rho ,t)+z_{\rho}(x,\rho ,t)=0\quad \mbox{\rm in}\ \Omega\times
(0,1)\times (0,+\infty ),\label{e2dzeta}\\
& &u (x,t) =0\quad \mbox{\rm on}\ \partial\Omega\times
(0,+\infty)\label{e3d}\\
& &\eta^t (x,s) =0\quad \mbox{\rm in}\ \partial\Omega\times
(0,+\infty), \ t\ge 0,\label{e4d}\\
& &z(x,0,t)=u_t(x,t)\quad \mbox{\rm in}\ \Omega\times
(0,+\infty), \label{e4dzeta}\\
& &u(x,0)=u_0(x)\quad \mbox{\rm and}\quad u_t(x,0)=u_1(x)\quad \hbox{\rm
in}\ \Omega,\label{e5d}\\
& & \eta^0(x,s)=\eta_0(x,s) \quad \mbox{\rm in}\ \partial\Omega\times
(0,+\infty), \label{e6d}\\
& &z(x,\rho , 0)=z^0(x,-\tau\rho) \quad x\in\Omega, \ \rho\in (0,1),\label{e7d}
\end{eqnarray}
where
\begin{equation}\label{datiinizd}
\begin{array}{l}
u_0(x)=u_0(x,0), \quad x\in\Omega,\\
u_1(x)=\frac {\partial u_0}{\partial t}(x,t)\vert_{t=0},\quad x\in\Omega,\\
\eta_0(x,s)=u_0(x,0)-u_0(x,-s),\quad x\in\Omega,\  s\in (0,+\infty),\\
z^0(x,s)=\frac {\partial u_0}{\partial t}(x,s),\quad x\in\Omega,\ s\in (-\tau ,0).
\end{array}
\end{equation}

Let us denote 
${\mathcal U}:= (u,u_t,\eta^t, z)^T.$ The we can rewrite    problem (\ref{e1d})--(\ref{e7d})
in the abstract form
\begin{equation}\label{abstractd}
\left\{
\begin{array}{l}
{\mathcal U}^\prime={\mathcal A} {\mathcal U},\\
{\mathcal U}(0)=(u_0,u_1, \eta_0, z^0)^T,
\end{array}
\right.
\end{equation}
where the operator ${\mathcal A}$ is defined by
\begin{equation}\label{Operator}
{\mathcal A}\left (
\begin{array}{l}
u\\v\\w\\z
\end{array}
\right )
:=\left (
\begin{array}{l}
v\\
(1-\tilde\mu)\Delta u+\int_0^{\infty}\mu (s)\Delta w(s)ds-k z(\cdot, 1)\\
-w_s+v
\\
-\tau^{-1} z_\rho
\end{array}
\right )\,,
\end{equation}
with domain

\begin{equation}\label{dominioOpd}
\begin{array}{l}
{\mathcal D}({\mathcal A}):=\left\{
\ (u,v,\eta, z)^T\in   H^1_0(\Omega)\times H^1_0(\Omega)
\times L^2_{\mu}((0,+\infty);H^1_0(\Omega))\times
 H^1((0,1);L^2(\Omega))\, :\right.\\\medskip
\hspace{3 cm}
v=z(\cdot, 0), \ (1-\tilde\mu)u+\int_0^\infty \mu (s)\eta (s) ds \in H^2(\Omega)\cap H^1_0(\Omega),\\
\hspace{9 cm}\left.
\ \eta_s\in  L^2_{\mu}((0,+\infty);H^1_0(\Omega))
\right\},
\end{array}
\end{equation}
where
$L^2_{\mu}((0, \infty);H^1_0(\Omega ))$ is the Hilbert space
of $H^1_0-$ valued functions on $(0,+\infty),$
endowed with the inner product
$$\langle \varphi, \psi\rangle_{L^2_{\mu}((0, \infty);H^1_0(\Omega ))}=
\int_{\Omega}\left (\int_0^\infty \mu (s)\nabla \varphi (x,s)\nabla \psi (x,s) ds\right )dx.
$$

\noindent
Denote by ${\mathcal H}$
the Hilbert space  

$${\mathcal H}= 
H^1_0(\Omega)\times L^2(\Omega)\times L^2_{\mu}((0, \infty);H^1_0(\Omega ))
\times L^2 ((0,1);L^2(\Omega)),$$
equipped
  with the inner product

\begin{equation}\label{innerd}
\begin{array}{l}
\left\langle
\left (
\begin{array}{l}
u\\
v\\
w\\
z
\end{array}
\right ),\left (
\begin{array}{l}
\tilde u\\
\tilde v\\
\tilde w\\
\tilde z
\end{array}
\right )
\right\rangle_{\mathcal H}
:= \displaystyle{
 (1-\tilde\mu )\int_\Omega \nabla u\nabla\tilde u dx + \int_\Omega v\tilde v dx +
\int_{\Omega} \int_0^\infty \mu (s)\nabla w\nabla\tilde w ds dx}\\
\medskip\hspace{8 cm}\displaystyle{
+\int_0^1\int_\Omega z(x,\rho) \tilde z(x,\rho)\, dxd\rho.}
\end{array}
\end{equation}

Combining the ideas from \cite{Pruss93} with the ones from  \cite{NPSicon06} (see also \cite{ANP10}), we can prove that the operator ${\mathcal A}$ generates a strongly continuous semigroup (${\mathcal A}-cI$ is dissipative for a sufficiently large constant $c>0$) and therefore
the next existence result holds.

\begin{Proposition}\label{WellP}
For any initial datum ${\mathcal U}_0\in {\mathcal H}$ there exists a unique solution
${\mathcal U}\in C([0,+\infty), {\mathcal H})$  of problem $(\ref{abstractd}).$
Moreover, if ${\mathcal U}_0\in {\mathcal D}({\mathcal A}),$ then
$${\mathcal U}\in C([0,+\infty), {\mathcal D}({\mathcal A}))\cap C^1([0,+\infty), {\mathcal H}).$$
\end{Proposition}

\bigskip

Let us define the energy $F$ of problem $(\ref{1.1d})-(\ref{1.3d})$ as
\begin{equation}\label{energyd}
\begin{array}{l}
\displaystyle{
F(t)=F(u,t):=\frac 1 2 \int_{\Omega} u_t^2(x,t) dx
+\frac {1-\tilde\mu} 2 \int_{\Omega}\vert \nabla u(x,t)\vert^2 dx}\\
\hspace{1 cm}\displaystyle{
+\frac 1 2 \int_0^{+\infty } \int_{\Omega}\mu(s)\vert \nabla\eta^t(s)\vert^2 ds dx+\frac {\theta\vert k\vert e^{\tau}} 2\int_{t-\tau}^t e^{-(t-s)}\int_{\Omega}u_t^2(x,s) ds dx,}
\end{array}
\end{equation}
where $\theta$ is any real constant satisfying
\begin{equation}\label{theta}
\theta >1.
\end{equation}

We will prove the following exponential stability result.

\begin{Theorem}\label{stab2}
For any $\theta >1$ in the definition $(\ref{energyd}),$
there exists a positive constant $k_0$ such that for  $k$ satisfying
$\vert k\vert <k_0$ there is $\sigma >0$ such that
\begin{equation}\label{exponentiald}
F(t)\le F(0) e^{1-\sigma t},\quad t\ge 0;
\end{equation}
for every solution of problem $(\ref{1.1d})-(\ref{1.3d}).$
The constant $k_0$ depends only on the kernel $\mu (\cdot)$ of the memory term, on the time delay $\tau$ and on the domain $\Omega.$
\end{Theorem}

To prove our stability result we will make use of the
following result
result of
 Pazy (Theorem 1.1 in Ch. 3 of \cite{pazy}).
\begin{Theorem}\label{Pazy}
Let $X$ be a Banach space and let $A$ be the infinitesimal generator of a $C_0$ semigroup $T(t)$ on $X,$ satisfying $\Vert T(t)\Vert \le Me^{\omega t}.$
If $B$ is a bounded linear operator on $X$ then $A+B$ is the infinitesimal generator of a $C_0$ semigroup $S(t)$ on $X,$ satisfying
$\Vert S(t)\Vert\le M e^{(\omega+M\Vert B\Vert )t}\,.$
\end{Theorem}

Moreover, we will use the
following lemma (see Th. 8.1 of \cite{Komornikbook}).

\begin{Lemma}\label{Vilmos}
Let $V(\cdot )$ be a non negative decreasing function defined on $[0,+\infty).$ If
$$
\int_S^{+\infty}V(t) dt \le C V(S)\,\quad\forall S>0\,,
$$
for some constants $C>0,$ then
$$
V(t)\le V(0) \exp \left (1-\frac{t}{C} \right ),\quad\forall\ t\ge 0\,.
$$
\end{Lemma}

\begin{Remark}\label{nodelay}
{\rm
Observe that the well--posedness result in the case $\tau =0,$ namely viscoelastic wave equation with standard frictional damping or anti--damping, directly follows from Theorem \ref{Pazy}.
Furthermore, from Theorem
\ref{Pazy} we can also deduce an exponential stability estimate under a suitable 
smallness assumption on  $\vert k\vert.$  
Indeed, for $\vert k\vert$ small, we can look at problem (\ref{1.1d})--(\ref{1.3d})(with $\tau =0$) as a perturbation of the wave equation with only the viscoelastic damping.
And it is by now well-known that for the last model an exponential decay estimate is available (see e.g. \cite{GiorgiRiveraPata}).
}
\end{Remark}

\section{Stability results}
\label{Stab1}\hspace{5mm}
\setcounter{equation}{0}

In this section we will prove Theorem \ref{stab2}.

\noindent
In order to study the stability properties of problem
(\ref{1.1d})--(\ref{1.3d}), we look at an auxiliary problem (cf. \cite{SCL12})
which is {\sl near} to this one and more easier to deal with.
Then, let us consider the system

\begin{eqnarray}
& &u_{tt}(x,t) -\Delta u (x,t)+
\int_0^\infty \mu (s)\Delta u(x,t-s) ds\nonumber \\
& &\hspace{3 cm}
 +\theta \vert k\vert e^\tau u_t(x,t)+
k u_t(x,t-\tau )=0\quad \mbox{\rm in}\ \Omega\times
(0,+\infty)\label{a.1d}\\
& &u (x,t) =0\quad \mbox{\rm on}\ \partial\Omega\times
(0,+\infty)\label{a.2d}\\
& &u(x,t)=u_0(x, t)\quad \hbox{\rm
in}\ \Omega\times (-\infty, 0]. \label{a.3d}
\end{eqnarray}

First of all we show that the energy, defined by (\ref{energyd}),
of any solution of the auxiliary problem is not increasing.

\begin{Proposition}\label{Fdecreasing}
For every solution of problem $(\ref{a.1d})-(\ref{a.3d})$ the energy $F(\cdot)$ is not increasing
and the following estimate holds
\begin{equation}\label{stimaF}
\begin{array}{l}
\displaystyle{
F^\prime(t)\le \frac 1 2 \int_0^\infty
\int_{\Omega}\mu^\prime(s)\vert \nabla \eta^t(x,s)\vert^2 dx ds
}\\\medskip
\hspace{1 cm}\displaystyle{
- \frac {\vert k\vert (\theta e^{\tau}-1)} 2 \int_{\Omega} u^2_t(x,t) dx - \frac {\vert k\vert (\theta -1)} 2\int_{\Omega}
u^2_t(x,t-\tau) dx}\\
\hspace{2 cm} \displaystyle{
 -\frac {\theta \vert k\vert e^{\tau}} 2\int_{t-\tau}^t e^{-(t-s)}\int_{\Omega}
u^2_t(x,s) dx ds\,.
}\\
\end{array}
\end{equation}
\end{Proposition}

\begin{Remark}{\rm
Note that the energy $F(\cdot)$ of solutions of the original problem
$(\ref{1.1d})-(\ref{1.3d})$
is not in general decreasing.
}
\end{Remark}

\noindent
{\bf Proof of Proposition \ref{Fdecreasing}.}
Differentiating (\ref{energyd}) we have
$$
\begin{array}{l}
\displaystyle{F^\prime (t)=\int_{\Omega}u_t(x,t)u_{tt}(x,t) dx +
(1-\tilde\mu )\int_{\Omega}\nabla u(x,t)\nabla u_t(x,t) dx}
\\
\hspace{1 cm} \displaystyle{
+\int_0^\infty \int_{\Omega } \mu (s)\nabla \eta^t(x,s)\nabla \eta^t_t(x,s) dx ds
+\frac{\theta\vert k\vert e^{\tau} } 2\int_{\Omega}u_t^2(x,t) dx
}\\
\hspace{1 cm} \displaystyle{
-\frac {\theta\vert k\vert} 2 \int_{\Omega} u^2_t(x,t-\tau) dx -\frac {\theta\vert k\vert e^{\tau}} 2
\int_{t-\tau}^t e^{-(t-s)}\int_{\Omega}u^2_t(x,s) dx ds\,.
}
\end{array}
$$
Then, integrating by parts and using
(\ref{e2d}) and the boundary condition (\ref{a.2d}),
$$
\begin{array}{l}
\displaystyle{F^\prime (t)=\int_{\Omega}u_t(x,t)[u_{tt}(x,t)
-(1-\tilde\mu )\Delta u(x,t)]dx}
\\
\hspace{1 cm} \displaystyle{
+\int_0^\infty \int_{\Omega } \mu (s)\nabla \eta^t(x,s)(\nabla u_t(x,t)-
\nabla
 \eta^t_s(x,s)) dx ds
+\frac{\theta\vert k\vert e^{\tau} } 2\int_{\Omega}u_t^2(x,t) dx
}\\
\hspace{1 cm} \displaystyle{
-\frac {\theta\vert k\vert} 2 \int_{\Omega} u^2_t(x,t-\tau) dx -\frac {\theta\vert k\vert e^{\tau}} 2
\int_{t-\tau}^t e^{-(t-s)}\int_{\Omega}u^2_t(x,s) dx ds\,.
}
\end{array}
$$

By using equations $(\ref{a.1d}), (\ref{a.2d}),$
after integration by parts, we deduce

$$
\begin{array}{l}
\displaystyle{F^\prime (t)=\int_{\Omega}u_t(t)\Big [
-\int_0^\infty \mu(s)\Delta u(x,t-s)+\tilde\mu \Delta u(x,t)}\\
\hspace{6.5 cm}\displaystyle{
-
\theta\vert k\vert e^{\tau}u_t(x,t)-ku_t(x,t-\tau)
\Big ]dx}
\\
\hspace{0.7 cm} \displaystyle{
+\int_0^\infty \int_{\Omega } \mu (s)\nabla \eta^t(x,s)\nabla u_t(x,t) dx ds
+\frac 1 2 \int_0^\infty\int_{\Omega}\mu^\prime (s)\vert\nabla \eta^t(x,s)\vert^2 dx ds}\\
\hspace{0.7 cm} \displaystyle{+\frac{\theta\vert k\vert e^{\tau} } 2\int_{\Omega}u_t^2(x,t) dx
-\frac {\theta\vert k\vert} 2 \int_{\Omega} u^2_t(x,t-\tau) dx -\frac {\theta\vert k\vert e^{\tau}} 2
\int_{t-\tau}^t e^{-(t-s)}\int_{\Omega}u^2_t(x,s) dx ds\,
}\\\hspace{0.5 cm}
=\displaystyle{
-\theta \vert k\vert e^{\tau}\int_{\Omega}u_t^2(x,t) dx-k\int_{\Omega}u_t(x,t)
u_t(x,t-\tau ) dx+\frac{\theta\vert k\vert e^{\tau} } 2\int_{\Omega}u_t^2(x,t) dx
}\\
\hspace{1 cm}\displaystyle{
-\frac {\theta\vert k\vert} 2 \int_{\Omega} u^2_t(x,t-\tau) dx
+\frac 1 2 \int_0^\infty\int_{\Omega}\mu^\prime (s)\vert\nabla \eta^t(x,s)\vert^2 dx ds}\\
\hspace{1 cm}\displaystyle{
 -\frac {\theta\vert k\vert e^{\tau}} 2
\int_{t-\tau}^t e^{-(t-s)}\int_{\Omega}u^2_t(x,s) dx ds\,.
}
\end{array}
$$

Now, using Cauchy-Schwarz inequality we obtain $(\ref{stimaF})$.\qed

\begin{Corollary}\label{coroFdecreasing}
For every solution of problem $(\ref{a.1d})-(\ref{a.3d})$, we have
\begin{equation}\label{corostimaF}
\displaystyle{
-\frac 1 2 \int_S^T\int_0^\infty
\int_{\Omega}\mu^\prime(s)\vert \nabla \eta^t(x,s)\vert^2 dx ds
\leq F(S),
}
\end{equation}
and then by  the condition $\mu'(t)\leq -\alpha \mu(t)$ we directly get
\begin{equation}\label{Al33}
\frac 12 \int_S^T\int_0^\infty \mu (s)\int_{\Omega} \vert \nabla \eta^t (x,s)
\vert^2 dx ds dt \le
\frac 1\alpha  F(S)\,.
\end{equation}
\end{Corollary}
\noindent{\bf Proof.}
As each term of the right-hand side of
(\ref{stimaF}) is non positive, we directly get that
$$-\frac 1 2 \int_S^T\int_0^\infty
\int_{\Omega}\mu^\prime(s)\vert \nabla \eta^t(x,s)\vert^2 dx ds
\leq  \int_S^T(-F^\prime (t)) dt\leq F(S).
\quad\quad
\qed $$

\begin{Theorem}\label{stimaint}
For any $\theta >1$ in the definition $(\ref{energyd}),$
there exist positive constants $C$
and $\overline{k},$ depending on $\mu,$ $\Omega$ and $\tau,$
such that if $\vert k\vert <\overline{k}$ then for any solution of problem $(\ref{a.1d})-(\ref{a.3d})$ the following estimate holds
\begin{equation}\label{integrale}
\int_S^{+\infty}F(t) dt \le C F(S)\,\quad\forall S>0\,.
\end{equation}
\end{Theorem}

In order to prove Theorem \ref{stimaint} we need some preliminary results.
Our proof relies in many points on \cite{ACS} but we have to perform all computations
because, in order to extend the exponential estimate related to the perturbed problem (\ref{a.1d})--(\ref{a.3d})
to the original problem (\ref{1.1d})--(\ref{1.3d}) we need to determine carefully all involved constants.
From the definition of the energy we deduce

\begin{equation}\label{Al1}
\begin{array}{l}
\displaystyle{\int_S^T
F(t)dt =\frac 1 2 \int_S^T\int_{\Omega} u_t^2(x,t) dx dt
+\frac {1-\tilde\mu } 2 \int_S^T\int_{\Omega}\vert \nabla u(x,t)\vert^2 dx dt}\\
\hspace{2.5 cm}\displaystyle{
+\frac 1 2\int_S^T \int_0^\infty \int_{\Omega}\mu (s)\vert \nabla \eta^t(x,s)\vert^2dx  ds dt}\\
\hspace{3 cm} \displaystyle{
+
\frac {\theta\vert k\vert e^{\tau}} 2\int_S^T\int_{t-\tau}^t e^{-(t-s)}\int_{\Omega}u_t^2(x,s) ds dx dt
\,.}
\end{array}
\end{equation}
Now, as in \cite{ACS} we will use multiplier arguments in order to bound the right--hand side of
(\ref{Al1}).
We note that we could not apply the same arguments directly to our original problem since the energy is not
decreasing.

In the following we will denote by $C_P$ the  Poincar\'e constant, namely the smallest positive constant such that
\begin{equation}\label{Poincare}
\int_\Omega w^2(x) dx \le C_P\int_\Omega \vert \nabla w(x)\vert^2 dx,\quad \forall\ w\in H^1_0(\Omega).
\end{equation}

\begin{Lemma}\label{lemma1}
Assume
\begin{equation}\label{prima}
\vert k\vert <\frac{1-\tilde \mu}  {2 C_P(\theta e^{\tau}+1)}.
\end{equation}
Then,
for any $T\ge S\ge 0$ we have
\begin{equation}\label{Al2}
(1-\tilde\mu )\int_S^T\int_{\Omega} \vert\nabla u(x,t)\vert^2 dx dt
\le C_0\int_S^T\int_{\Omega}u^2_t(x,t) dx dt+ C_1 F(S),
\end{equation}
with
\begin{equation}\label{C0eC1}
C_0=2+\theta\vert k\vert e^{\tau}\,,
\quad\quad
C_1=4\Big (1+
\frac {\tilde \mu}{\alpha (1-\tilde\mu)} +\frac {C_P}{1-\tilde\mu }+ \frac 1 {2
(\theta -1)}
\Big )\,.
\end{equation}
\end{Lemma}

\noindent{\bf Proof.} Multiplying equation (\ref{a.1d}) by $u$ and integrating on $\Omega\times [S,T]$ we have
$$
\begin{array}{l}
\displaystyle{
\int_S^T\int_{\Omega} [u_{tt}(x,t)-\Delta u (x,t)+\int_0^\infty\mu (s)
\Delta u (x,t-s) ds}\\
\hspace{4 cm}\displaystyle{
+\theta \vert k\vert e^{\tau}
u_t(x,t)+ku_t(x,t-\tau)] u(x,t) dx dt =0}\,.
\end{array}
$$

So, integrating by parts and using the boundary condition (\ref{a.2d}), we get
$$
\begin{array}{l}
\displaystyle{
-\int_S^T\int_{\Omega}u^2_t(x,t) dx dt+\int_S^T\int_{\Omega}\vert\nabla u (x,t)\vert^2 dx dt+\Big [\int_{\Omega}u (x,t)u_t(x,t)dx  \Big ]_S^T
}\\
\hspace{2 cm}
\displaystyle{+\theta\vert k\vert e^{\tau}\int_S^T\int_{\Omega}u (x,t)u_t(x,t) dx dt
+k\int_S^T\int_{\Omega}u (x,t)u_t(x,t-\tau ) dx dt}\\
\hspace{2 cm} \displaystyle{
-\tilde \mu\int_S^T\int_{\Omega}\vert\nabla u (x,t)\vert^2 dx dt
+\int_S^T\int_{\Omega}\int_0^\infty\mu (s)\nabla u (x,t)\nabla\eta^t(x,s)ds dx dt=0\,,
}
\end{array}
$$
where we used (\ref{eta}).

Then,
\begin{equation}\label{Al3}
\begin{array}{l}
\displaystyle{
(1-\tilde\mu )
\int_S^T\int_{\Omega}\vert\nabla u (x,t)\vert^2 dx dt}\\
\hspace{1.5 cm}
\displaystyle{= \int_S^T \int_{\Omega}u^2_t(x,t) dx dt-\Big [\int_{\Omega}u (x,t)
u_t(x,t) dx  \Big ]_S^T}\\
\hspace{2 cm}\displaystyle{
-\theta\vert k\vert e^{\tau}\int_S^T\int_{\Omega}u (x,t)u_t(x,t) dx dt
-k\int_S^T\int_{\Omega}u (x,t)u_t(x,t-\tau ) dx dt
}\\
\hspace{2.5 cm} \displaystyle{
-\int_S^T\int_{\Omega}\int_0^\infty\mu (s)\nabla u (x,t)
\nabla\eta^t(x,s)ds dx dt\,.
}
\end{array}
\end{equation}
In order to estimate the integral
$$
\int_S^T\left\vert \int_{\Omega}\int_0^\infty\mu (s)\nabla \eta^t(x,s)\nabla u(x,t)ds dx\right\vert  dt\,,$$
we note that, for all $\varepsilon >0,$
\begin{equation}\label{Al4}
\begin{array}{l}
\displaystyle{
\int_S^T\Big ( \int_{\Omega} \vert \nabla u (x,t)\vert^2 dx \Big )^{1/2}
\int_0^\infty\mu (s)\Big ( \int_{\Omega} \vert \nabla\eta^t(x,s)\vert^2 dx \Big )^{1/2} ds dt
}\\
\displaystyle{\le
\frac {\varepsilon} 2 \int_S^T\int_{\Omega} \vert \nabla u (x,t)\vert^2 dx dt+\frac 1 {2\varepsilon} \int_S^T\Big [\int_0^\infty\mu (s)\Big (
\int_{\Omega} \vert \nabla\eta^t(x,s)\vert^2 dx \Big )^{1/2} ds \Big ]^2 dt\,.
}
\end{array}
\end{equation}
We have
$$
\begin{array}{l}
\displaystyle{
\int_S^T\Big [ \int_0^\infty\mu (s)\Big ( \int_{\Omega} \vert \nabla \eta^t(x,s)\vert^2 dx \Big )^{1/2} ds \Big ]^2 dt
}\\
\hspace{1 cm}
\displaystyle{
\le\int_S^T\Big (\int_0^\infty \mu (s) ds\Big )\Big (\int_0^\infty\mu (s)\int_{\Omega}\vert \nabla \eta^t(x,s)\vert^2 dx ds\Big ) dt
}\\
\hspace{1 cm}
\displaystyle{=\tilde\mu\int_S^T\int_0^\infty\mu (s)\int_{\Omega}\vert\nabla\eta^t(x,s)\vert^2 dx ds dt\,.}
\end{array}
$$

\noindent
Therefore, recalling the estimate (\ref{Al33}), we obtain
\begin{equation}\label{Al5}
\int_S^T\Big [ \int_0^\infty\mu (s)\Big ( \int_{\Omega} \vert \nabla \eta^t(x,s)\vert^2 dx \Big )^{1/2} ds \Big ]^2 dt
\le \frac {2\tilde\mu}
{\alpha }F(S)\,.
\end{equation}

\noindent
Then, (\ref{Al4}) and (\ref{Al5}) give
\begin{equation}\label{Al6}
\begin{array}{l}
\displaystyle{
\int_S^T\left\vert \int_{\Omega}\int_0^\infty\mu (s)(\nabla u (x,t-s)-
\nabla u(x,t))\cdot\nabla
u (x,t) ds dx\right\vert  dt}\\
\hspace{2 cm}\le \displaystyle{
\frac {\varepsilon}{2}\int_S^T\int_{\Omega} \vert\nabla u (x,t)\vert^2 dx dt+\frac {\tilde \mu} {\alpha\varepsilon} F(S)\,.
}
\end{array}
\end{equation}
Now observe that
\begin{equation}\label{Al7}
F(t)\ge \frac 1 2 \int_{\Omega} u_t^2(x,t) dx +\frac {1-\tilde\mu } 2
\int_{\Omega}\vert \nabla u(x,t)\vert^2 dx\,.
\end{equation}
Then, from (\ref{Al7}),
\begin{equation}\label{Al8}
\frac 1 2\int_{\Omega} \vert \nabla u (x,t)\vert^2 dx \le \frac {F(t)}{1-\tilde
\mu}\,,
\end{equation}
and also, from Poincar\'e's inequality,
\begin{equation}\label{Al9}
\frac 1 2 \int_{\Omega} \vert u (x,t)\vert^2 dx \le \frac {C_P} 2\int_{\Omega} \vert \nabla u (x,t)\vert^2 dx \le\frac {C_P} {1-\tilde\mu } F(t)\,.
\end{equation}
Using the above inequalities
\begin{equation}\label{Al10}
\left\vert \int_{\Omega} u_t(x,t)u(x,t) dx\right\vert
\le \frac 1 2 \int_{\Omega}u_t^2(x,t) dx +\frac 1 2 \int_{\Omega}u^2(x,t)dx\le F(t)\Big (1+\frac {C_P}{1-\tilde\mu }\Big )\,.
\end{equation}
Therefore,
\begin{equation}\label{Al11}
-\Big [\int_{\Omega}u_t(x,t)u (x,t) dx \Big ]_S^T\le 2 F(S)\Big (
1+\frac {C_P}{1-\tilde\mu }\Big )\,,
\end{equation}
where we used also the fact that $F$ is decreasing.
Using (\ref{Al6}), (\ref{Al11}) and Cauchy--Schwarz's inequality in order to bound the terms in the right--hand side of (\ref{Al3}) we have that for any $\varepsilon >0,$
$$
\begin{array}{l}
\displaystyle{
(1-\tilde\mu )\int_S^T \int_{\Omega} \vert\nabla u (x,t)\vert^2 dx dt\le
\int_S^T\int_{\Omega}u^2_t(x,t) dx dt+
\frac {\varepsilon} 2\int_S^T\int_{\Omega}\vert\nabla u (x,t)\vert^2 dx dt
}\\
\hspace{0.4 cm} +\displaystyle{
\frac {\tilde\mu } {\alpha \varepsilon}F(S)+
2\Big (
1+\frac {C_P}{1-\tilde\mu }\Big )F(S)+\frac {\theta\vert k\vert e^{\tau}} 2 \int_S^T\int_{\Omega}u^2 (x,t)dx dt
}\\
\hspace{0.4 cm} +\displaystyle{
\frac  {\theta\vert k\vert e^{\tau}} 2 \int_S^T\int_{\Omega}u^2_t (x,t)dx dt
+\frac {\vert k\vert} 2 \int_S^T\int_{\Omega}u^2 (x,t)dx dt+
\frac {\vert k\vert} 2 \int_S^T\int_{\Omega}u^2_t(x,t-\tau)dx dt\,.
}
\end{array}
$$
Therefore, from Poincar\'e's inequality,

$$
\begin{array}{l}
\displaystyle{
(1-\tilde\mu )\int_S^T \int_{\Omega} \vert\nabla u (x,t)\vert^2 dx dt\le \Big (1+\frac  {\theta\vert k\vert e^{\tau}}
 2  \Big ) \int_S^T\int_{\Omega}u^2_t(x,t) dx dt}\\
\hspace{1 cm}\displaystyle{+
\frac {\varepsilon+(\theta e^{\tau}+1)\vert k\vert C_P} 2\int_S^T\int_{\Omega}\vert\nabla
u (x,t)\vert^2 dx dt+\frac {\tilde\mu } {\alpha\varepsilon}F(S)
}\\
\hspace{1 cm}\displaystyle{+
2\Big (
1+\frac {C_P}{1-\tilde\mu }\Big )F(S)+\frac {\vert k\vert} 2 \int_S^T\int_{\Omega}u^2_t(x,t-\tau)dx dt\,.
}
\end{array}
$$

\noindent
Now, observe that from (\ref{stimaF}),
\begin{equation}\label{Al12}
\begin{array}{l}
\displaystyle{
\frac{\vert k\vert } 2 \int_S^T \int_{\Omega} u^2_t(x,t-\tau) dx dt
=\frac {1} {\theta -1}
\frac {\vert k\vert(\theta -1) } 2 \int_S^T\int_\Omega u^2_t(x,t-\tau) dx dt
}\\
\hspace{1 cm}
\displaystyle{\le \frac {1} {\theta -1}\int_S^T(-F^\prime(t)) dt\leq\frac 1 {\theta -1 } F(S)\,.
}
\end{array}
\end{equation}
Now, choose $\varepsilon= \frac{1-\tilde\mu } 2.$ Thus, using (\ref{prima}) and  also (\ref{Al12})
we obtain
$$
\begin{array}{l}
\displaystyle{
(1-\tilde\mu )\int_S^T \int_{\Omega} \vert\nabla u (x,t)\vert^2 dx dt\le 2\Big (1+\frac {\theta \vert k\vert e^{\tau}} 2   \Big )\int_S^T\int_{\Omega} u_t^2(x,t) dx dt
}\\
\hspace{1 cm}\displaystyle{
+4\Big (1+\frac{\tilde\mu }{\alpha (1-\tilde\mu )}
+\frac{C_P}{1-\tilde\mu }
+\frac{1} {2(\theta -1)}
\Big ) F(S)\,,
}
\end{array}
$$
that is (\ref{Al2}) with constants $C_0, C_1$ given by (\ref{C0eC1}).\qed

\begin{Lemma}\label{lemma2}
For any $T\ge S\ge 0,$ the following identity holds:
\begin{equation}\label{Al14}
\begin{array}{l}
\displaystyle{
\tilde\mu\int_S^T\int_{\Omega} u^2_t(x,t) dx dt
=\Big [
\int_{\Omega}u_t(x,t)\int_0^\infty\mu (s)\eta^t(x,s) ds dx
\Big ]_S^T
}\\
\hspace{1 cm}-\displaystyle{
\int_S^T\int_{\Omega}u_t(x,t)\int_0^\infty\mu^\prime (s)\eta^t(x,s) ds dx dt
}\\
\hspace{1 cm}+\displaystyle{ (1-\tilde\mu )
\int_S^T \int_{\Omega}\nabla u(x,t)
\int_0^\infty \mu (s) \nabla\eta^t(x,s) ds dx dt
}\\
\hspace{1 cm}+\displaystyle{\int_S^T\int_{\Omega}\Big \vert
\int_0^\infty\mu (s) \nabla\eta^t (x,s)ds
\Big\vert^2 dx dt
}\\
\hspace{1 cm}+\displaystyle{
\theta\vert k\vert e^{\tau}\int_S^T\int_{\Omega} u_t(x,t)\int_0^\infty \mu (s)
\eta^t(x,s) ds dx dt
}\\
\hspace{1 cm}+\displaystyle{
k\int_S^T\int_{\Omega} u_t(x,t-\tau)\int_0^\infty \mu (s) \eta^t(x,s) ds dx dt\,.
}
\end{array}
\end{equation}
\end{Lemma}

\noindent {\bf Proof.} We multiply equation (\ref{a.1d}) by
$\int_0^\infty\mu (s)\eta^t(x,s) ds$ and integrate by parts on
$[S,T]\times \Omega$. We obtain

\begin{equation}\label{Al15}
\begin{array}{l}
\displaystyle{
\int_S^T\int_{\Omega} \Big\{u_{tt}(x,t)-\Delta u(x,t)+\int_0^\infty\mu (s)
\Delta u(x,t-s) ds +k u_t(x,t-\tau)+\theta\vert k\vert e^{\tau } u_t(x,t)\Big\}
}\\
\hspace{3 cm}\displaystyle{\times
\Big\{
\int_0^\infty \mu (s)\eta^t(x,s) ds
\Big\} dx dt =0\,.
}
\end{array}
\end{equation}
Integrating by parts, we have
\begin{equation}\label{Al16}
\begin{array}{l}
\displaystyle{
\int_S^T\int_{\Omega} u_{tt}(x,t)\int_0^\infty\mu (s)
\eta^t(x,s) ds dx dt
}\\
\hspace{1 cm}
\displaystyle{=\Big [
\int_{\Omega}u_t (x,t)\int_0^\infty\mu (s)\eta^t(x,s) ds dx
\Big ]_S^T
}\\
\hspace{1.3 cm}
\displaystyle{
-\int_S^T\int_{\Omega} u_t(x,t)\int_0^\infty\mu (s)(u_t(x,t)-\eta^t_s (x,s)) ds dx dt
}\\
\hspace{1 cm}
\displaystyle{=\Big [
\int_{\Omega}u_t (x,t)\int_0^\infty\mu (s)\eta^t(x,s) ds dx
\Big ]_S^T
}\\
\hspace{1.3 cm}
\displaystyle{
-\tilde\mu\int_S^T\int_{\Omega}u_t^2(x,t) dx dt
-\int_S^T\int_{\Omega} u_t(x,t)\int_0^\infty\mu^\prime (s)\eta^t (x,s) ds dx dt
\,.
}
\end{array}
\end{equation}
Moreover,
\begin{equation}\label{Al17}
\begin{array}{l}
\displaystyle{
\int_S^T\int_{\Omega} \Big (
-\Delta u (x,t)+\int_0^\infty\mu (s)\Delta u (x,t-s) ds
\Big )\int_0^\infty\mu (s) \eta^t(x,s) ds dx dt
}\\
\hspace{1 cm}
\displaystyle{
=\int_S^T\int_{\Omega}\nabla u (x,t) \int_0^\infty\mu (s)\nabla \eta^t (x,s)ds dx dt
}\\
\hspace{1.3 cm}
\displaystyle{-\int_S^T\int_{\Omega}\int_0^\infty\mu (s)\nabla u (x,t-s) ds
 \int_0^\infty\mu (s)\nabla \eta^t (x,s) ds dx dt
}\\
\hspace{1 cm}
\displaystyle{
=\int_S^T\int_{\Omega}\nabla u (x,t) \int_0^\infty\mu (s)\nabla \eta^t (x,s) ds dx dt
}\\
\hspace{1.3 cm}
\displaystyle{+\int_S^T\int_{\Omega}\int_0^\infty\mu (s)(\nabla u (x,t)
-\nabla u (x,t-s))
 ds
 \int_0^\infty\mu (s)\nabla \eta^t(x,s)  ds dx dt
}\\
\hspace{1.3 cm}
\displaystyle{
-\tilde\mu\int_S^T\int_{\Omega}\nabla u(x,t)
\int_0^\infty\mu (s)
\nabla \eta^t (x,s)  ds dx dt
}\\
\hspace{1 cm}
\displaystyle{
=(1-\tilde\mu \int_S^T
\int_{\Omega}\nabla u (x,t) \int_0^\infty\mu (s)\nabla \eta^t (x,s) ds dx dt
}\\
\hspace{1.3 cm}
\displaystyle{
+\int_S^T\int_{\Omega}\Big \vert
\int_0^\infty\mu (s)\nabla \eta^t (x,s) ds
\Big\vert^2 dx dt\,.
}
\end{array}
\end{equation}
Using (\ref{Al16}) and (\ref{Al17}) in (\ref{Al15}) we obtain (\ref{Al14}).\qed
\begin{Lemma}\label{lemma3}
Assume
\begin{equation}\label{Al26}
\vert k\vert <\frac{\tilde\mu }{2\theta}e^{-\tau}\,.
\end{equation}
Then,
for any $T\ge S>0$ and for any $\varepsilon >0$ we have
\begin{equation}\label{Al18}
\int_S^T\int_{\Omega}u_t^2(x,t) dx dt\le \varepsilon  \int_S^T\int_{\Omega}
\vert\nabla u (x,t)\vert^2 dx dt+ C_2 F(S)\,,
\end{equation}
where
the constant  $C_2:=C_2(\varepsilon)$ is defined by
\begin{equation}\label{C2}
C_2=
\frac 4 {\tilde\mu}\Big (
1+\frac 1 2 \frac 1 {\theta -1}+\frac{\mu (0)}{\tilde\mu }C_P\Big )
+4C_P
+\frac {2}{\alpha } \Big (2+\frac {(1-\tilde\mu )^2} {\tilde\mu\varepsilon}
+ C_P\vert k\vert (\theta e^{\tau}+1)\Big )
\,.
\end{equation}
\end{Lemma}
\noindent {\bf Proof.} In order to prove Lemma \ref{lemma3} we have to estimate the terms of the right-hand side of (\ref{Al14}).
First we have,
$$
\begin{array}{l}
\displaystyle{
\Big \vert
\int_{\Omega}u_t(x,t)\int_0^\infty\mu (s) \eta^t(x,s) ds dx\Big\vert
}\\
\hspace{1 cm}\displaystyle{\le
\int_0^\infty\mu (s)\Big (
\int_{\Omega}\vert u_t(x,t)\vert \vert \eta^t(x,s)\vert dx
\Big )ds
}\\
\hspace{1 cm}\displaystyle{
\le\int_0^\infty\mu (s)\Big (\int_{\Omega} u_t^2(x,t) dx\Big )^{1/2}
\Big (
\int_{\Omega}(\eta^t (x,s))^2 dx
\Big )^{1/2}
ds
}\\
\hspace{1 cm}\displaystyle{
\le\frac 1 2\int_{\Omega} u^2_t(x,t) dx +\frac 1 2 \Big (
\int_0^\infty\mu (s)\Big (\int_{\Omega}(\eta^t(x,s))^2
dx\Big )^{1/2} ds
\Big  )^2\,.
}
\end{array}
$$
Then, recalling (\ref{energyd}) and using H\"{o}lder's inequality, we deduce
\begin{equation}\label{Al19}
\begin{array}{l}
\displaystyle{
\Big \vert
\int_{\Omega}u_t(x,t)\int_0^\infty\mu (s) \eta^t (x,s) ds dx\Big\vert
}\\
\hspace{0.3 cm}
\displaystyle{
\le F(t)+\frac {C_P} 2 \Big (
\int_0^\infty\mu (s) \Big (\int_{\Omega}\vert \nabla\eta^t (x,s)\vert^2 dx
\Big )^{1/2} ds \Big )^2
}\\
\hspace{0.3 cm}
\displaystyle{\le F(t)+\frac {C_P} 2 \tilde\mu
\int_0^\infty\mu (s) \int_{\Omega}\vert \nabla\eta^t (x,s)\vert^2 dx
ds\le F(t)(1+C_P\tilde\mu )\,. }
\end{array}
\end{equation}
Therefore,
\begin{equation}\label{Al20}
\Big [\int_{\Omega} u_t(x,t)\int_0^\infty\mu (s)\eta^t (x,s)ds dx \Big ]_S^T\le 2 (1+C_P\tilde\mu ) F(S)\,.
\end{equation}
Now we proceed to estimate the second term in the right--hand side of (\ref{Al14}). For any $\delta >0$ we have
$$
\begin{array}{l}
\displaystyle{\left\vert
\int_S^T\int_{\Omega}u_t(x,t)\int_0^\infty\mu^\prime (s)\eta^t (x,s)ds dx dt\right\vert
}\\
\hspace{1 cm}\displaystyle{\le
\int_S^T\Big (\int_{\Omega}u_t^2(x,t) dx \Big )^{1/2}
\Big ( \int_{\Omega}
\Big ( \int_0^\infty\mu^\prime (s) \eta^t (x,s) ds\Big )^2 dx
\Big )^{1/2} dt
}\\
\hspace{1 cm}\displaystyle{\le
\frac {\delta} 2 \int_S^T\int_{\Omega}u^2_t(x,t) dx dt
+\frac 1 {2\delta} \int_S^T
\int_{\Omega}
\Big ( \int_0^\infty\mu^\prime (s) \eta^t (x,s) ds\Big )^2 dx dt}\\
\hspace{1 cm}\displaystyle{\le
\frac {\delta} 2 \int_S^T\int_{\Omega}u^2_t(x,t) dx dt
}\\
\hspace{3 cm}\displaystyle{
+\frac 1 {2\delta} \int_S^T
\int_{\Omega}
\int_0^\infty (-\mu^\prime (s) )ds
 \int_0^\infty\vert \mu^\prime (s)\vert (\eta^t (x,s))^2 ds dx dt}\,,
\end{array}
$$
and then by Corollary \ref{coroFdecreasing}
\begin{equation}\label{Al21}
\begin{array}{l}
\displaystyle{\left\vert
\int_S^T\int_{\Omega}u_t(x,t)\int_0^\infty\mu^\prime (s)\eta^t (x,s) ds dx dt\right\vert
}\\
\hspace{1 cm} \displaystyle{
\le \frac \delta 2 \int_S^T\int_{\Omega}u_t^2(x,t) dx dt- \frac {\mu (0)}
{2\delta }C_P\int_S^T\int_0^\infty \mu^\prime (s)\int_{\Omega} \vert \nabla\eta^t (x,s)\vert^2 dx ds dt
}\\
\hspace{1 cm} \displaystyle{
\le \frac \delta 2 \int_S^T\int_{\Omega}u_t^2(x,t) dx dt+ \frac {\mu (0)}
{\delta }C_P F(S)\,.}
\end{array}
\end{equation}
Moreover, by (\ref{Al33}) we  have
\begin{equation}\label{Al23}
\begin{array}{l}
\displaystyle{
\int_S^T\int_{\Omega}\Big \vert \int_0^\infty\mu (s) \nabla \eta^t (x,s) ds \Big \vert^2 dx dt
}\\
\hspace{1 cm}
\displaystyle{\le
\int_S^T\int_{\Omega}\tilde\mu \int_0^\infty\mu (s) \vert \nabla \eta^t (x,s)\vert ^2ds  dx dt
}\\
\hspace{1 cm}\displaystyle{ \leq \frac {2\tilde \mu }
{\alpha }F(S)\,.
}
\end{array}
\end{equation}
Then, it results also
\begin{equation}\label{Al22}
\begin{array}{l}
\displaystyle{
\int_S^T \int_{\Omega} \nabla u (x,t)
\int_0^\infty\mu (s) \nabla\eta^t(x,s)ds dx dt
}\\
\hspace{1 cm}\displaystyle{
\le\int_S^T\Big\vert\int_{\Omega}\nabla u (x,t)\int_0^\infty\mu (s)
\nabla \eta^t (x,s) ds dx \Big\vert dt
}\\
\hspace{1 cm}
\displaystyle{
\le \frac {\varepsilon} 2\int_S^T\int_{\Omega}\vert\nabla u (x,t)\vert^2 dx dt
+\frac {\tilde\mu } {\alpha \varepsilon} F(S)\,.
}
\end{array}
\end{equation}

\noindent
Now we estimate the last two integrals in the right--hand side of (\ref{Al14}).
$$
\begin{array}{l}
\displaystyle{
\theta\vert k\vert e^{\tau}\int_S^T\int_{\Omega} u_t(x,t)\int_0^\infty\mu (s)
\eta^t (x,s) ds dx dt\hspace{4 cm}
}\\
\displaystyle{
+k\int_S^T\int_{\Omega} u_t(x,t-\tau)\int_0^\infty\mu (s)\eta^t (x,s) ds dx dt\hspace{1.7 cm}
}\\
\hspace{1 cm}
\displaystyle{\le
\frac {\vert k\vert} 2 \int_S^T\int_{\Omega}u^2_t(x,t-\tau ) dx dt
+\frac {\theta\vert k\vert e^{\tau}} 2 \int_S^T\int_{\Omega} u^2_t(x,t) dx dt}\\
\hspace{1.5 cm}\displaystyle{
+\frac {\vert k\vert(1+\theta e^{\tau})} 2\int_S^T\int_{\Omega} \Big (
\int_0^\infty\mu (s)
\eta^t (x,s) ds \Big )^2 dx dt
}\\
\hspace{1 cm}
\displaystyle{\le
\frac {\vert k\vert } 2 \int_S^T\int_{\Omega}u^2_t(x,t-\tau ) dx dt
+\frac {\theta\vert k\vert e^{\tau}}  2 \int_S^T\int_{\Omega} u^2_t(x,t) dx dt}\\
\hspace{1.5 cm}\displaystyle{+
\frac {\vert k\vert(1+\theta e^{\tau})} 2 C_P\tilde\mu \int_S^T\int_{\Omega} \int_0^\infty\mu (s)
\vert \nabla \eta^t (x,s)\vert^2 ds  dx dt
}\,.
\end{array}
$$
Therefore, recalling (\ref{Al33}) and (\ref{Al12}), we have
\begin{equation}\label{Al24}
\begin{array}{l}
\displaystyle{
\theta\vert k\vert e^{\tau}\int_S^T\int_{\Omega} u_t(x,t)\int_0^\infty\mu (s)
\eta^t (x,s) ds dx dt\hspace{4 cm}
}\\
\displaystyle{
+k\int_S^T\int_{\Omega} u_t(x,t-\tau)\int_0^\infty\mu (s)\eta^t(x,s) ds dx dt\hspace{1.7 cm}
}\\
\hspace{0.2 cm}
\displaystyle{
\le
\frac {\vert k\vert} 2 \int_S^T\int_{\Omega}u^2_t(x,t-\tau ) dx dt
+\frac {\theta\vert k\vert e^{\tau}} 2 \int_S^T\int_{\Omega} u^2_t(x,t) dx dt+C_P
(\vert k\vert (\theta e^{\tau}+1))\frac {\tilde\mu}{\alpha}F(S)
}\\
\hspace{0.2 cm}\displaystyle{
\le \frac {1} {\theta -1 } F(S)+\frac {\theta\vert k\vert e^{\tau}} 2
\int_S^T\int_{\Omega} u_t^2 (x,t) dx dt+
\frac {C_P\tilde\mu} {\alpha} \vert k\vert (\theta e^{\tau}+1 ) F(S)\,.
}
\end{array}
\end{equation}
Using (\ref{Al20})--(\ref{Al24}) in (\ref{Al14}) we obtain
\begin{equation}\label{Al25}
\begin{array}{l}
\displaystyle{
\Big ( \tilde\mu -\frac {\theta\vert k\vert e^{\tau}} 2-\frac\delta 2\Big )\int_S^T\int_{\Omega}
u^2_t(x,t) dx\le \frac {\varepsilon} 2(1-\tilde\mu )\int_S^T\int_{\Omega}
\vert\nabla u (x,t)\vert^2 dx dt}\\
\hspace{ 1 cm} \displaystyle {+
\frac {1} {\theta -1} F(S)+
2(1+C_P\tilde\mu ) F(S)+\frac{\mu (0)} \delta C_P F(S)
}\\
\hspace{1 cm}
\displaystyle{
+\frac {\tilde \mu } \alpha (\frac {1- \tilde\mu} {\varepsilon} +2 ) F(S)
+ C_P\tilde\mu \frac {\vert k\vert (\theta e^{\tau}+1)}{\alpha} F(S)
\,.
}
\end{array}
\end{equation}

Now, fix $\delta =\frac {\tilde\mu } 2 \,.$
Then, from (\ref{Al26}), for any $T\ge S>0,$ we have

\begin{equation}\label{Al28}
\begin{array}{l}
\displaystyle{
\int_S^T\int_{\Omega}u_t^2(x,t) dx dt\le
\frac\varepsilon {\tilde\mu} (1-\tilde\mu )\int_S^T\int_{\Omega}
\vert\nabla u(x,t)\vert^2 dx dt
}\\
\hspace{1 cm} \displaystyle{
+
\frac 2 {\tilde\mu}
\Big ( 2 (1+C_P\tilde\mu )+
\frac {1}{\theta -1} + 2\frac {\mu (0)}{\tilde\mu } C_P+
\frac{\tilde\mu}{\alpha } (2+ \frac {1-\tilde\mu }\varepsilon
+C_P\vert k\vert(\theta e^{\tau}+1)  )
\Big )F(S)},
\end{array}
\end{equation}
that is (\ref{Al18}) with constant $C_2$ as in (\ref{C2}).\qed
\begin{Lemma}\label{lemma4}
Assume
\begin{equation}\label{Al29}
\vert k\vert <\min\ \Big \{
\frac {1-\tilde\mu}{2C_P(\theta e^{\tau}+1)}, \frac {\tilde \mu}{2\theta}e^{-\tau}
\Big \}\,.
\end{equation}
Then, for any $T\ge S>0,$
\begin{equation}\label{Al30}
\frac {1-\tilde\mu } 2\int_S^T\int_{\Omega}\vert\nabla u (x,t)\vert^2 dx dt +\frac 12\int_S^T\int_{\Omega}u_t^2(x,t) dx dt
\le C^* F(S)\,,
\end{equation}
with
 \begin{equation}\label{Cstar}
C^*= C_0C_2+ C_1+C_2\,,
\end{equation}
where $C_0$ and  $C_1$ are the constants defined by $(\ref{C0eC1})$ and
\begin{equation}\label{C2bis}
\begin{array}{l}
\displaystyle{
C_2:= C_2\Big (\frac{1-\tilde\mu}{2(C_0+1)}  \Big )=\frac 4 {\tilde\mu}\Big (
1+\frac 1 2 \frac 1 {\theta -1}+\frac{\mu (0)}{\tilde\mu }C_P\Big )
+4C_P}\\\medskip
\hspace{4.5 cm}\displaystyle{
+\frac {2}{\alpha } \Big (2+(6+2\theta\vert k\vert e^{\tau})\frac {1-\tilde\mu } {\tilde\mu}
+ C_P\vert k\vert (\theta e^{\tau}+1)\Big )}
\,.
\end{array}
\end{equation}
\end{Lemma}

\noindent
{\bf Proof.} The assumptions of previous lemmas are verified. Thus,
we can use  (\ref{Al18}) in (\ref{Al2}). Then,
\begin{equation}\label{Al31}
\begin{array}{l}
\displaystyle{
(1-\tilde\mu )\int_S^T\int_{\Omega}\vert\nabla u (x,t)\vert^2 dx}\\
\hspace{1.5 cm}\displaystyle{
\le C_0\varepsilon \int_S^T\int_{\Omega}\vert\nabla u (x,t)\vert^2dx dt+ (C_0C_2+C_1) F(S)\,.}
\end{array}
\end{equation}
Therefore, from (\ref{Al18}) and (\ref{Al31}), we obtain
\begin{equation}\label{Al32}
\begin{array}{l}
\displaystyle{
\frac {1-\tilde\mu }2\int_S^T \int_{\Omega}\vert\nabla u (x,t)\vert^2 dx
+\frac 12\int_S^T\int_{\Omega}u^2_t(x,t) dx dt
}\\
\displaystyle{
\hspace{1 cm}\le
\frac \varepsilon 2 (C_0+1)\int_S^T\int_{\Omega} \vert\nabla u (x,t)\vert^2 dx dt +\frac 12 (C_0C_2+C_1+C_2) F(S)\,.
}
\end{array}
\end{equation}
Now, fix
$$\varepsilon =\frac {1-\tilde\mu } {2(C_0+1)}\,.$$
Then, from (\ref{Al32}) we deduce
$$
\begin{array}{l}
\displaystyle{
\frac {1-\tilde\mu } 4\int_S^T \int_{\Omega}\vert\nabla u (x,t)\vert^2 dx
+\frac 12\int_S^T\int_{\Omega}u^2_t(x,t) dx dt
}\\
\displaystyle{
\hspace{1 cm}\le
\frac 12 (C_0C_2+C_1+C_2) F(S)\,,
}
\end{array}
$$
where, from (\ref{C2}) with the above choice of $\varepsilon,$
$C_2$ is as in (\ref{C2bis}).
This clearly implies (\ref{Al30}) with $C^*$ as in (\ref{Cstar}).\qed

\noindent {\bf Proof of Theorem \ref{stimaint}.}
Notice also that (\ref{corostimaF})  directly implies that
\begin{equation}\label{Al34}
\frac {\theta\vert k\vert e^{\tau}} 2 \int_S^T\int_{t-\tau}^t e^{-(t-s)}\int_{\Omega} u^2_t(x,s) dx ds dt\le -\int_S^T F^\prime (t) dt \le F(S)\,.
\end{equation}
Let us define $\overline{k}$ as
\begin{equation}\label{kbarrato}
\overline{k}:=
\min\ \Big \{
\frac {1-\tilde\mu}{2C_P(\theta e^{\tau}+1)}, \frac {\tilde \mu}{2\theta}e^{-\tau}
\Big \}\,.
\end{equation}

\noindent
Then, if $\vert k\vert <\overline{k},$ using  (\ref{Al30}), (\ref{Al33}) and (\ref{Al34}) in (\ref{Al1}), we obtain
$$
\int_S^T F(t) dt\le C^* F(S)+\frac 1\alpha F(S)+ F(S)\,.
$$
Therefore (\ref{integrale}) is verified with
\begin{equation}\label{C}
C=C^*+1+\frac 1\alpha\,,
\end{equation}
where $C^*$ is as in (\ref{Cstar}) with $C_0, C_1$ and $C_2$
defined in (\ref{C0eC1}) and (\ref{C2bis}).\qed

{\bf Proof of Theorem \ref{stab2}}
From Theorem \ref{stimaint} and Lemma \ref{Vilmos},
it follows that for any solution of the auxiliary problem
$(\ref{a.1d})-(\ref{a.3d})$ if $\vert k\vert <\overline{k},$
we have
\begin{equation}\label{expaux}
F(t)\le F(0)e^{1-\tilde\sigma t},\quad t\ge 0,
\end{equation}
 with
\begin{equation}\label{sigmatilde}
\tilde\sigma:=\frac 1 C,
\end{equation}
where $C$ is as in (\ref{C}).

From this and Theorem \ref{Pazy} we deduce that Theorem \ref{stab2} holds,
with $\sigma:=\tilde\sigma -e\theta \vert k\vert e^{\tau},$
if
$$
-\tilde\sigma +e\theta\vert k\vert e^\tau <0,
$$
that is if the delay parameter $k$ satisfies
\begin{equation}\label{finale}
\vert k\vert <g (\vert k\vert ):=\frac 1 {Ce\theta e^{\tau}},
\end{equation}
with $C:=C(\vert k\vert )$ defined in (\ref{C}).
Now observe that (\ref{finale}) is satisfied for $k=0$ because $g(0)>0.$
Moreover, by recalling the definitions of the constants $C_0, C_1, C_2$ and $C^*,$ used to define $C,$
we note that $g:[0,\infty)\rightarrow (0,\infty )$ is a continuous decreasing function  satisfying
$$g(\vert k\vert)\rightarrow 0
  \quad \mbox{ \rm for}\quad \vert k\vert \rightarrow\infty.$$
Thus, there exists a unique constant $\hat k >0$ such that $\hat k= g(\hat k).$
We can then conclude that for any $\theta$ in the definition (\ref{energyd})
of the energy $F(\cdot),$
inequality (\ref{finale}) is satisfied for every $k$ with
$$\vert k\vert < k_0=\min\{ \hat k, \overline{k}\}.\hspace{2 cm}\qed$$

\begin{Remark} \label{explicitk0}
{\rm
We can compute an explicit lower bound  for  $k_0.$
Indeed (\ref{finale}) may be rewritten as
$$
\vert k\vert \theta e^{\tau +1}\left (C^*+1+\frac 1 {\alpha}\right )<1.$$
Then, from (\ref{Cstar}), we have
\begin{equation}\label{q1}
[1+1/\alpha +C_2(C_0+1)+C_1]\theta e^{\tau +1}\vert k\vert <1,
\end{equation}
that is
\begin{equation}\label{q2}
\begin{array}{l}
\displaystyle{h(\vert k\vert):=
\left\{1+\frac 1 {\alpha}+\left [
\frac 4 {\tilde\mu}\Big (
1+\frac 1 2 \frac 1 {\theta -1}+\frac{\mu (0)}{\tilde\mu }C_P\Big )
+4C_P \right.\right.}\\
\displaystyle{
\hspace{1 cm}\left.
+\frac {2}{\alpha } \Big (2+(6+2\theta\vert k\vert e^{\tau})\frac {1-\tilde\mu } {\tilde\mu}
+ C_P\vert k\vert (\theta e^{\tau}+1)\Big )\right ](3+\theta \vert k\vert
e^{\tau})}
\\\displaystyle{ \left. \hspace{2.5 cm} +
4\Big (1+
\frac {\tilde \mu}{\alpha (1-\tilde\mu)} +\frac {C_P}{1-\tilde\mu }+ \frac 1 {2
(\theta -1)}
\Big )
\right\}\theta e^{\tau+1}\vert k\vert } <1.
\end{array}
\end{equation}
Now, we use the assumption $\vert k\vert <\overline {k}$
with $\overline{k}$ defined in (\ref{kbarrato}) in order to majorize the left--hand side of (\ref{q2}), $h(\vert k\vert ),$ with a linear function.
We have
\begin{equation}\label{q3}
\begin{array}{l}
\displaystyle{h(\vert k\vert)\le
\left\{1+\frac 1 {\alpha}+\left [
\frac 4 {\tilde\mu}\Big (
1+\frac 1 2 \frac 1 {\theta -1}+\frac{\mu (0)}{\tilde\mu }C_P\Big )
+4C_P \right.\right.}\\
\displaystyle{
\hspace{1 cm}\left.
+\frac {2}{\alpha } \Big (2+(6+\tilde\mu )\frac {1-\tilde\mu } {\tilde\mu}
+ \frac {1-\tilde\mu }2\Big )\right ](3+\tilde\mu /2)}
\\\displaystyle{ \left. \hspace{2.5 cm} +
4\Big (1+
\frac {\tilde \mu}{\alpha (1-\tilde\mu)} +\frac {C_P}{1-\tilde\mu }+ \frac 1 {2
(\theta -1)}
\Big )
\right\}\theta e^{\tau+1}\vert k\vert },
\end{array}
\end{equation}
from which follows
$$h(\vert k\vert )\le \left (1+\frac 1 {\alpha }\gamma_1+\gamma_2\right )\theta \vert k\vert e^{\tau+1},$$
with
\begin{eqnarray*}
\gamma_1&=&\gamma_1(\tilde\mu )= 4\frac {\tilde \mu} {1-\tilde\mu }-8+\frac{36}{\tilde\mu }-\frac{23} 2\tilde\mu -\frac 3 2 \tilde\mu^2,
\\
\gamma_2&=&\gamma_2(\mu(0),\tilde \mu, \theta , C_P)\\
&=&
\displaystyle{
6+12C_P +\frac 3 {\theta -1}+\frac {12} {\tilde\mu}+\frac 6 {\tilde\mu (\theta -1)} }
\\
&&\hspace{0.5cm}+\displaystyle{12\frac {\mu (0)}{\tilde\mu^2}C_P+2\frac {\mu (0)}{\tilde\mu}C_P
+2C_P\tilde\mu  +\frac {4C_P}{1-\tilde\mu }.}
\end{eqnarray*}
Then, we deduce the following explicit lower bound
\begin{equation}\label{explicit}
k_0\geq \frac {e^{-(\tau+1)}} {\theta (1+\frac 1 {\alpha}\gamma_1+\gamma_2)},
\end{equation}
with $\gamma_1,\gamma_2$ as before.
For example, if we take
$$\mu(t)=e^{-2t},$$
then $\tilde\mu =1/2$ and so, fixing $\theta =2,$ we can compute
$\gamma_1=\frac{495} 8,\quad \gamma_2=45+ 73 C_P.$
Hence,
for this particular choice of the memory kernel, we obtain
$$k_0\geq \frac{8 e^{-(\tau +1)}}{1231+1168 C_P}.$$
}
\end{Remark}

\begin{Remark}\label{nodelay2}
{\rm
In the case $\tau=0$ and  $k<0,$
namely viscoelastic wave equation with anti--damping, we can simplify previous arguments. Indeed, the absence of time delay allows us to take $\theta =1$
obtaining an exponential stability estimate under the condition
$$\vert k\vert <  
 \left (C_1+3C_2+\frac 1 {\alpha }\right )^{-1}\frac 1 e,$$
where
$$
C_1=4\Big (1+
\frac {\tilde \mu}{\alpha (1-\tilde\mu)} +\frac {C_P}{1-\tilde\mu }
\Big )$$
and
$$
C_2=\frac 2 {\tilde\mu}\Big (
2+\frac{\mu (0)}{\tilde\mu }C_P\Big )
+4C_P
+\frac {2}{\alpha } \Big (2+6\frac {1-\tilde\mu } {\tilde\mu}
\Big ).$$
}
\end{Remark}


\begin{thebibliography}{10}

\bibitem{ACS}
F.~Alabau-Boussouira, P.~Cannarsa and D.~Sforza.
\newblock Decay estimates for second order evolution equations with memory.
\newblock {\em J. Funct. Anal.}, 254:1342--1372, 2008.

\bibitem{AC}
F.~Alabau-Boussouira, P.~Cannarsa.
\newblock A new method for proving sharp energy decay rates for memory-dissipative evolution equations for a quasi-optimal class of kernels. 
\newblock {\em C. R. Acad. Sci. Paris, S\'er. I}, 347:867--872, 2009.


\bibitem{ANP10}
K.~Ammari, S.~Nicaise and C.~Pignotti.
\newblock Feedback boundary stabilization of wave
equations with interior delay.
\newblock {\em Systems Control Lett.}, 59:623--628, 2010.


\bibitem{Dafermos}
C.~M.~Dafermos.
\newblock Asymptotic stability in viscoelasticity.
\newblock {\em Arch. Rational Mech. Anal.}, 37:297--308, 1970.

\bibitem{Datko}
R.~Datko.
\newblock Not all feedback stabilized hyperbolic systems are robust
with respect to small time delays in their feedbacks.
\newblock {\em SIAM J. Control Optim.}, 26:697--713, 1988.

\bibitem{DLP}
R.~Datko, J.~Lagnese and M.~P.~Polis.
\newblock An example on the effect of time delays in boundary feedback stabilization of
wave equations.
\newblock  {\em SIAM J. Control Optim.}, 24:152--156, 1986.

\bibitem{FZ}
P.~Freitas and E.~Zuazua.
\newblock Stability results for the wave equation with indefinite damping.
\newblock {\em J. Differential Equations}, 132:338--352, 1996.

\bibitem{GiorgiRiveraPata}
C.~Giorgi, J.~E.~Mu\~{n}oz Rivera and V.~Pata.
\newblock Global Attractors for a Semilinear Hyperbolic Equation in Viscoelasticity.
\newblock {\em  J. Math. Anal. Appl.}, 260:83--99, 2001.


\bibitem{Guesmia}
A.~Guesmia.
\newblock Well--posedness and exponential stability of an abstract evolution
equation with infinite memory and time delay.
\newblock {\em IMA J. Math. Control Inform.}, 30:507--526, 2013.







\bibitem{Komornikbook}
V.~Komornik.
\newblock {\em Exact controllability and stabilization, the multiplier method},
  volume~36 of {\em RMA}.
\newblock Masson, Paris, 1994.













\bibitem{Rivera}
J.E.~ Mun\~{o}z Rivera and A.~ Peres Salvatierra.
\newblock
Asymptotic behaviour of the energy in partially viscoelastic materials.
\newblock {\em Quart. Appl. Math.}, 59:557--578, 2001.



\bibitem{NPSicon06}
S.~Nicaise and C.~Pignotti.
\newblock Stability and instability results of the wave equation with a delay term in the boundary
or internal feedbacks.
\newblock {\em SIAM J. Control Optim.}, 45:1561--1585, 2006.



\bibitem{NiPi2013}
S.~Nicaise and C.~Pignotti.
\newblock Stabilization of second--order evolution
       equations with time delay.
\newblock {\em Math. Control Signals Syst.}, DOI 10.1007/s00498-014-0130-1, 2014.



\bibitem{Pata}
V.~Pata.
\newblock
Exponential stability in linear viscoelasticity with almost flat memory kernels.
\newblock {\em Commun. Pure Appl. Anal.}, 9:721--730, 2010.
\bibitem{pazy}
A.~Pazy.
\newblock {\em Semigroups of linear operators and applications to partial differential equations}, Vol. 44 of {\em Applied Math. Sciences.} Springer-Verlag, New York, 1983.

\bibitem{SCL12}
C.~Pignotti.
\newblock A note on stabilization of locally damped wave equations with time delay.
\newblock {\em Systems and Control Lett.}, 61:92--97. 2012. 

\bibitem{Pruss93}  J. Pr\"{u}ss.
\newblock {\em  Evolutionary Integral Equations and Applications}, Monogr. Math., vol. 87, Birkh\"a user-Verlag, Basel, 1993.



\bibitem{XYL}
G.~Q.~Xu, S.~P.~Yung and L.~K.~Li.
\newblock Stabilization of wave systems with input delay in the boundary control.
\newblock {\em ESAIM Control Optim. Calc. Var.}, 12(4):770--785,  2006.



\end{thebibliography}
 \end{document}